\documentclass[twoside,11pt]{article}
\title{} \author{} \date{}
\usepackage{amssymb,latexsym,times}
\newtheorem{te}{Theorem}[section]
\newtheorem{fac}[te]{Fact}

\newtheorem{lem}[te]{Lemma}

\newtheorem{ex}[te]{Example}

\def\dok{\noindent{\bf Proof. }}
\def\kdok{\hfill $\Box$ \par \vspace*{2mm} }
\def\upar{\!\uparrow}
\def\downar{\!\downarrow}
\def\Lim{\lim\nolimits}
\newcommand{\Conv}{\mathop{\rm Conv}\nolimits}
\newcommand{\Br}{\mathop{{\mathrm {Br}}}}
\newcommand{\ro}{\mathop{{\mathrm {ro}}}}
\begin{document}
\thispagestyle{plain}
\begin{center}
           {\large \bf
           \uppercase{
                      A CONVERGENCE ON BOOLEAN ALGEBRAS GENERALIZING\\[1mm]
                         THE CONVERGENCE ON THE ALEKSANDROV CUBE
                     }
                     }
\end{center}

\begin{center}
{\small \bf Milo\v s S.\ Kurili\'c and Aleksandar Pavlovi\'c}\\[2mm]
{\small  Department of Mathematics and Informatics, University of Novi Sad,\\
Trg Dositeja Obradovi\'ca 4, 21000 Novi Sad, Serbia \\
e-mail: milos@dmi.uns.ac.rs, apavlovic@dmi.uns.ac.rs}
\end{center}
\begin{abstract}
\noindent
We compare the forcing related properties of a complete Boolean algebra ${\mathbb B}$ with the properties
of the convergences $\lambda _{\mathrm s}$ (the algebraic convergence) and $\lambda_{\mathrm{ls}}$ on ${\mathbb B}$
generalizing the convergence on the Cantor and Aleksandrov cube respectively. In particular we show that
$\lambda_{\mathrm{ls}}$ is a topological convergence iff forcing by ${\mathbb B}$ does not produce new reals and
that $\lambda_{\mathrm{ls}}$ is weakly topological if ${\mathbb B}$ satisfies condition $(\hbar )$
(implied by the ${\mathfrak t}$-cc).
On the other hand, if $\lambda_{\mathrm{ls}}$ is a weakly topological convergence, then ${\mathbb B}$ is
a $2^{\mathfrak h}$-cc algebra or in some generic extension the distributivity number of the ground model
is greater than or equal to the tower number of the extension.
So, the statement ``The convergence $\lambda_{\mathrm{ls}}$ on the collapsing algebra ${\mathbb B}=\ro ({}^{<\omega }\omega _2)$
is weakly topological" is independent of ZFC.


\end{abstract}

\section{Introduction}

The object of our study is the interplay between the forcing related properties of a complete Boolean algebra
${\mathbb B}$ and the properties of convergence structures defined on ${\mathbb B}$.
In Section 3 we observe the algebraic convergence $\lambda _{\mathrm s}$, generalizing the convergence on the Cantor cube
and generating the sequential topology ${\mathcal O}_{\lambda_{\mathrm{s} }}$ introduced by Maharam
and investigated in the context of the von Neumann and Maharam's Problem.
In the rest of the paper we investigate the convergence $\lambda_{\mathrm{ls} }$, introduced in Section 4 as
a natural generalization of
the convergence on the Aleksandrov cube.

Concerning the context of our research, first we note that the
topology ${\mathcal O}_{\lambda_{\mathrm{ls} }}$ (generated by the convergence $\lambda_{\mathrm{ls} }$) and its
dual ${\mathcal O}_{\lambda_{\mathrm{li} }}$ generate the sequential topology ${\mathcal O}_{\lambda_{\mathrm{s} }}$,
for the algebras ${\mathbb B}$ belonging to a wide class including Maharam algebras \cite{KuPaDEBR}.
Second, we mention some related results. If ${\mathbb B}$ is a complete Boolean algebra, let
the convergences $\lambda _i : {\mathbb B}^\omega \rightarrow P({\mathbb B})$, for $i\in \{0,1,2,3,4\}$,
be defined by
$$
\lambda_{i}(x)=\left\{
\begin{array}{cl}
\{b_{4}(x)\} &\mbox { if } b_{i}(x)=b_{4}(x),\\
\emptyset &\mbox { if } b_{i}(x)<b_{4}(x),
\end{array}
\right.
$$
where $x=\langle x_n \rangle\in {\mathbb B}^\omega$, $\tau
_x= \{ \langle \check{n}, x_n \rangle : n\in \omega \}$ is the
corresponding ${\mathbb B}$-name for a subset of $\omega$ and
\begin{eqnarray*}
b_{0}(x)&=& \|\tau_x \mbox { is cofinite} \| = \liminf x, \\
b_{1}(x)&=& \|\tau_x \mbox { is old infinite} \| , \\
b_{2}(x)&= &\|\tau_x \mbox { contains an old infinite subset} \|, \\
b_{3}(x)&=& \|\tau_x \mbox { is infinite and non-splitting} \|,\\
b_{4}(x)&=& \|\tau_x \mbox { is infinite} \| = \limsup x .
\end{eqnarray*}
Then, by \cite{KuPaNSJOM} and \cite{KuPaDEBR}, $\lambda _{\mathrm s}=\lambda _0$ and
$\lambda_{\mathrm{ls} }= \bar\lambda _2 = \bar\lambda _3 = \bar\lambda _4$, where $\bar \lambda$ is the closure of
a convergence $\lambda $ under (L2).
Also $\lambda _1 \leq \lambda _2 \leq\lambda _3 \leq\lambda _4 $ and these four convergences
generate the same topology, ${\mathcal O}_{\lambda_{\mathrm{ls} }}$. So
we have the following diagram ($\lambda ^*$ denotes the closure of
a convergence $\lambda $ under (L3), see Section 2).

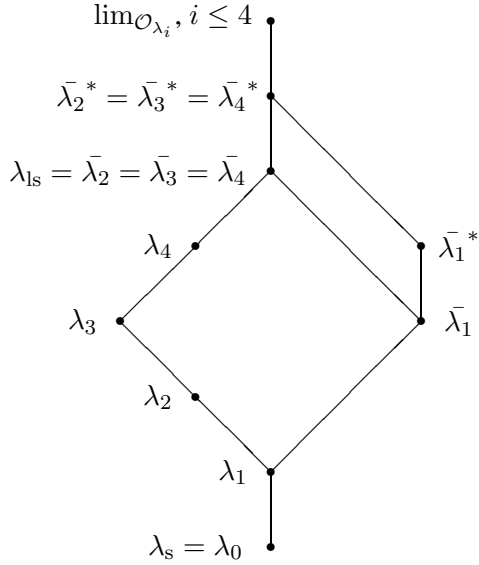
\unitlength 1mm 
\linethickness{0.4pt}
\ifx\plotpoint\undefined\newsavebox{\plotpoint}\fi 
\begin{picture}(130,75)(0,0)

\put(60,0){\line(0,1){10}}
\put(40,30){\line(1,-1){20}}
\put(60,10){\line(1,1){20}}
\put(40,30){\line(1,1){20}}
\put(60,50){\line(1,-1){20}}
\put(80,30){\line(0,1){10}}
\put(60,60){\line(1,-1){20}}
\put(60,50){\line(0,1){20}}
\put(60,0){\circle*{1}}
\put(60,10){\circle*{1}}
\put(50,20){\circle*{1}}
\put(40,30){\circle*{1}}
\put(80,30){\circle*{1}}
\put(80,40){\circle*{1}}
\put(50,40){\circle*{1}}
\put(60,50){\circle*{1}}
\put(60,60){\circle*{1}}
\put(60,70){\circle*{1}}


\put(50,0){\makebox(0,0)[cc]{$\lambda _{{\mathrm{s}}}=\lambda _0$}}
\put(55,10){\makebox(0,0)[cc]{$\lambda _1$}}
\put(45,20){\makebox(0,0)[cc]{$\lambda _2$}}
\put(35,30){\makebox(0,0)[cc]{$\lambda _3$}}
\put(85,30){\makebox(0,0)[cc]{$\bar{\lambda _1}$}}
\put(85,40){\makebox(0,0)[cc]{$\bar{\lambda _1}^*$}}
\put(45,40){\makebox(0,0)[cc]{$\lambda _4$}}
\put(41,50){\makebox(0,0)[cc]{$\lambda _{{\mathrm{ls}}}=\bar{\lambda _2}=\bar{\lambda _3}=\bar{\lambda _4}$}}
\put(45,60){\makebox(0,0)[cc]{$\bar{\lambda _2}^* = \bar{\lambda _3}^* = \bar{\lambda _4}^* $}}
\put(47,70){\makebox(0,0)[cc]{$\lim _{{\mathcal O}_{\lambda _i}}$, $i\leq 4$}}

\end{picture}
\vspace{2mm}

\noindent
Now we mention some related results from \cite{KuPaNSJOM} and \cite{KuPaDEBR}.
The property that ${\mathbb B}$ does not produce
new reals by forcing is equivalent to each of the following conditions: $\lambda _1 = \lambda _2 $,
$\lambda _1 = \lambda _4 $, $\lambda _2 = \lambda _3 $, $\lambda _2 = \lambda _4 $,
$\bar{\lambda} _1 = \lambda _{\mathrm{ls} } $, $\bar{\lambda} _1$ is a topological convergence.
The property that ${\mathbb B}$ does not produce splitting reals is equivalent to the equality $\lambda _3 = \lambda _4$,
which holds if the convergence $\bar{\lambda} _1$ is weakly topological.

Our notation is mainly standard. So,  $\omega$ denotes the set of
natural numbers, $Y^X$ the set of all functions $f: X \rightarrow Y$
and $\omega^{\uparrow \omega}$ the set of all strictly increasing
functions from $\omega$ into $\omega$.
By $|X|$ we denote the cardinality of the set $X$
and, if $\kappa$ is a cardinal, then $[X]^\kappa=\{A\subset
X:|A|=\kappa\}$  and $[X]^{<\kappa}=\{A\subset X : |A|<\kappa\}$.
By $\mathfrak{c}$ we denote the cardinality of the continuum.
For subsets $A$ and $B$ of $\omega$ we write $A \subset^*B$ if $A
\setminus B$ is a finite set and $A\varsubsetneq^*B$ denotes $A
\subset^*B$ and $B\not \subset^* A$. The set $S$ {\bf splits} the
set $A$ if the sets $A\cap S$ and $A\setminus S$ are infinite.
${\mathcal S}\subset [\omega]^\omega$ is called a {\bf splitting
family} if each set $A\in [\omega]^\omega$ is split by some element
of ${\mathcal S}$ and $\mathfrak{s}$ is the minimal size of a
splitting family (the {\bf splitting number}). A set $P$ is a {\bf
pseudointersection} of a family ${\mathcal T}\subset[\omega]^\omega$
if   $P\subset ^*T$ for each $T \in {\mathcal T}$. A family
${\mathcal T}\subset [\omega]^\omega$ is a {\bf tower} if $\langle
{\mathcal T},^*\!\!\varsupsetneq\rangle $ is a well-ordered set and
${\mathcal T}$ has no pseudointersection. The {\bf tower number},
$\mathfrak{t}$, is the minimal size of a tower in $[\omega]^\omega$.
If $\langle {\mathbb P},\leq \rangle $ is a partial order, a subset
$D\subset {\mathbb P}$ is called {\bf dense} if $\forall p \in
{\mathbb P}$ $\exists d \in D$  $d \leq p$ and $D$ is called {\bf
open} if $p\leq q \in D$ implies $p \in D$. The {\bf distributivity
number}, $\mathfrak{h}$, is the minimal size of a family of dense
open subsets of the order $\langle [\omega]^\omega,\subset \rangle $
whose intersection is not dense. More information on invariants of the continuum
the reader can find in \cite{Douw84}.

If ${\mathbb B}$ is a Boolean algebra and  $A\subset{\mathbb B}$
let $A\upar=\{ b \in \mathbb{B} : \exists a \in A \; a\leq b\}$;
instead of $\{b\}\upar$ we will write $b \upar$. Clearly, $A\upar =\bigcup_{a \in A} a\upar$
and we will say that a set $A$ is {\bf upward closed} iff $A=A\upar$.
In a similar way we define $A\downar$, $b\downar$ and {\bf downward closed} sets.

\section{Topological preliminaries}

A {\bf sequence} in a set $X$
is each function $x:\omega \rightarrow X$; instead of $x(n)$ we
usually write $x_n$ and also $x=\langle x_n:n\in \omega\rangle $.
The {\bf constant sequence} $\langle a,a,a,\dots\rangle$ is denoted
by $\langle a \rangle$.  If $f \in \omega^{\uparrow \omega}$, the
sequence $y =x \circ f$  is said to be a {\bf subsequence} of the
sequence $x$ and we write $y \prec x$.

If $\langle X, {\mathcal O}\rangle $ is a topological space, a point
$a\in X$ is said to be a {\bf limit point of a sequence} $x\in
X^{\omega }$ (we will write: $x \rightarrow_{\mathcal O}a$) iff each
neighborhood $U$ of $a$ contains all but finitely many members of
the sequence. A space $\langle X,{\mathcal O}\rangle $ is called
{\bf sequential} iff a set $A \subset X$ is closed whenever it
contains each limit of each sequence in $A$.

If $X$ is a non-empty set, each mapping $\lambda : X^{\omega
}\rightarrow P(X)$ is a {\bf convergence} on $X$ and the mapping
$u_{\lambda }:P(X)\rightarrow P(X)$, defined by $u_{\lambda }(A)=
\bigcup_{x\in A^{\omega }}\lambda (x)$, the {\bf operator of
sequential closure} determined by $\lambda$. A convergence $\lambda$
satisfying $|\lambda(x)|\leq 1$, for each sequence $x$ in $X$, is called a
{\bf Hausdorff convergence}. If $\lambda _1$ is another convergence
on $X$, then we will write $\lambda \leq \lambda _1$ iff $\lambda
(x) \subset \lambda _1 (x)$, for each sequence $x\in X^\omega$.
Clearly, $\leq$ is a partial ordering on the set $\Conv (X) = \{
\lambda :\lambda \mbox{ is a convergence on } X \}$.

If $\langle X,{\mathcal O}\rangle $ is a topological space, then the
mapping $\lim_{\mathcal O}:X^\omega\rightarrow P(X)$ defined by
$\lim_{\mathcal O}(x)=\{ a\in X :x \rightarrow_{\mathcal O}a \}$ is
{\bf the convergence on} $X$ {\bf determined by the topology}
${\mathcal O}$ and for the operator $\lambda=\lim_{\mathcal O}$ we have (see \cite{Eng85})\\[-2mm]

(L1) $\forall a \in X\;\;  a \in \lambda (\langle a\rangle )$;

(L2) $\forall x \in X^\omega\; \forall y \prec x\; \lambda (x) \subset \lambda (y)$;

(L3) $\forall x \in X^\omega \;
     \forall a \in X\;
     ((\forall y \prec x\;
     \exists z\prec y \;
     a \in \lambda(z)) \Rightarrow a \in \lambda(x))$.\\[-3mm]

\noindent
We will use the following facts which mainly belong to the topological folklore. Their proofs can be found in
\cite{KuPaNSJOM}.

\begin{fac}\label{T1200} \rm
If ${\mathcal O}_1$ and ${\mathcal O}_2 $ are topologies on a set $X$, then

(a) ${\mathcal O}_1\subset {\mathcal O}_2$ implies $\lim _{{\mathcal O}_2}\leq \lim _{{\mathcal O}_1}$.

(b) If ${\mathcal O}_1$ and ${\mathcal O}_2$ are sequential topologies and
     $\lim _{{\mathcal O}_1} = \lim _{{\mathcal O}_2}$, then ${\mathcal O}_1 ={\mathcal O}_2$.
\end{fac}
\noindent
A convergence $\lambda : X^{\omega }\rightarrow P(X)$ is called a
{\bf topological convergence} iff there is a topology ${\mathcal O}$
on $X$ such that $\lambda =\lim_{\mathcal O}$. The following fact
 shows that each convergence has
a minimal topological extension  and connects topological and
convergence structures.
\begin{fac}\rm\label{T1201}
Let $\lambda : X^{\omega }\rightarrow P(X)$ be a convergence on a
non-empty set $X$. Then

(a) There is the maximal topology ${\mathcal O}_{\lambda }$  on $X$
satisfying $\lambda  \leq \lim\nolimits_{{\mathcal O}_\lambda}$;

(b) $ {\mathcal O}_\lambda=
    \{ O\subset X : \forall x \in X^\omega \, (O \cap \lambda(x) \neq \emptyset
    \Rightarrow \exists n_0 \in \omega~ \forall n \geq n_0~ x_n \in O)\}$;

(c) $\langle X,{\mathcal O}_\lambda \rangle$ is a sequential space;

(d) ${\mathcal O}_\lambda = \{ X\setminus F : F\subset X \wedge
u_{\lambda }(F) =F \}$,
    if $\lambda$ satisfies (L1) and (L2);

(e) $\lim _{{\mathcal O}_{\lambda }} =\min \{ \lambda ' \in \Conv
(X) :
                                      \lambda ' \mbox{ is topological and } \lambda \leq \lambda ' \}$;

(f) ${\mathcal O}_{\lim _{{\mathcal O}_{\lambda }}} ={\mathcal
O}_{\lambda }$;

(g) If $\lambda _1 : X^{\omega }\rightarrow P(X)$ and $\lambda _1
\leq \lambda $,
    then ${\mathcal O} _{\lambda } \subset {\mathcal O} _{\lambda _1}$;

(h) $\lambda $ is a topological convergence iff $\lambda =\lim
_{{\mathcal O}_{\lambda }}$.
\end{fac}
If a convergence $\lambda$ satisfies conditions (L1) and (L2),
then the minimal closure of $\lambda$ under (L3) is
described in the following statement.
\begin{fac}\rm\label{T1247}
Let  $\lambda : X^{\omega }\rightarrow P(X)$ be a convergence satisfying (L1) and (L2). Then
the mapping $\lambda^* :X^\omega \rightarrow P(X)$ given by
$ \textstyle \lambda^*(y)
=\bigcap_{f\in \omega^{\uparrow \omega}} \bigcup_{g\in \omega^{\uparrow \omega}}{\lambda}(y \circ f \circ g)$
is the minimal convergence bigger than ${\lambda }$ and satisfying
(L1) - (L3). Hence $\lambda ^* \leq \lim _{{\mathcal O}_{\lambda }}$.
\end{fac}
A convergence $\lambda:X^\omega
\rightarrow P(X)$ will be called {\bf weakly-topological} iff it
satisfies conditions (L1) and (L2) and its (L3)-closure, $\lambda ^*$, is a topological convergence.
\begin{fac}\label{T1220} \rm
Let $\lambda:X^\omega \rightarrow P(X)$ be a convergence satisfying
(L1) and (L2).

(a) $\lambda$ is weakly topological iff $\lambda^*=\lim_{{\mathcal O}_\lambda }$,
that is iff for each $x\in X^\omega$ and $a\in X$ we have:
$
\textstyle a \in \lim_{{\mathcal O}_\lambda}(x) \Leftrightarrow
\forall y \prec x \;\; \exists z \prec y \;\; a \in \lambda(z)$;

(b) If $\lambda$ is a Hausdorff convergence, then $\lambda^*$ is
also a Hausdorff convergence and $\lambda^*=\lim_{{\mathcal
O}_\lambda}$, that is $\lambda$ is a weakly-topological convergence.
\end{fac}
\begin{fac}\label{T1206}\rm
Let $\lambda : X^{\omega }\rightarrow P(X)$ be  a convergence
satisfying (L1) and (L2) and let the mappings $u_{\lambda}^{\alpha
}: P(X) \rightarrow P(X)$, $\alpha \leq \omega _1$, be defined by
recursion in the following way: for $A\subset X$

$u_{\lambda}^0(A)=A$,

$u_{\lambda}^{\alpha +1}(A) = u_{\lambda }(u_{\lambda}^{\alpha
}(A))$ and

$u_{\lambda}^{\gamma }(A) = \bigcup _{\alpha <\gamma
}u_{\lambda}^{\alpha }(A)$, for a limit $\gamma \leq \omega _1$.

\noindent Then $u_{\lambda}^{\omega _1}$ is the closure operator in
the space $\langle X, {\mathcal O}_\lambda \rangle $.
\end{fac}

\section{The Cantor cube and the algebraic convergence}

First we recall that if $X_n$, $n\in \omega$, is a sequence of sets, then
$\liminf _{n\in \omega}X_n = \bigcup _{k\in \omega}\bigcap _{n\geq k}X_n
                           =  \{ x: x\in X_n \mbox{ for all but finitely many }n \}$
and
$\limsup _{n\in \omega}X_n  = \bigcap _{k\in \omega}\bigcup _{n\geq k}X_n
                          =  \{ x: x\in X_n \mbox{ for infinitely many }n \}$.
Clearly we have
\begin{fac} \rm \label{T1244}
Let $X_n$, $n\in \omega$, be a sequence of sets. Then

(a) $\liminf _{n\in \omega}X_n \subset \limsup _{n\in \omega}X_n$;

(b) If $X_n=X$, for each $n\geq k$, then $\liminf _{n\in \omega}X_n $ $= \limsup _{n\in \omega}X_n = X$.
\end{fac}
We remind the reader that if $\kappa $ is an infinite cardinal, then the Cantor cube of weight $\kappa $,
denoted by $\langle 2^\kappa, \tau _{\mathrm{C}}\rangle $, is the Tychonov product of $\kappa$ many copies
of the two point discrete space $2=\{ 0,1 \}$. We will identify the set $2^\kappa$ with the power set $P(\kappa )$
using the bijection $f : 2^\kappa\rightarrow P(\kappa )$ defined by $f(x)=x^{-1}[ \{ 1\} ]$.
\begin{fac} \rm \label{T1245}
Let $\langle x_n \rangle$ be a sequence in $2^\kappa$ and $x\in 2^\kappa$.
Then the following conditions are equivalent:

(a) $\langle x_n \rangle \rightarrow _{\tau _{\mathrm{C}}} x$,

(b) $\forall \alpha \in \kappa \; \exists k\in \omega \; \forall n\geq k \; x_n (\alpha )=x(\alpha )$,

(c) $\textstyle \liminf _{n\in \omega}X_n = \limsup _{n\in \omega}X_n =X$, where $X_n = f(x_n)$ and $X = f(x)$.

\noindent
The Cantor cube $\langle 2^\kappa, \tau_{\mathrm C}\rangle $ is a sequential space iff $\kappa=\omega$.
\end{fac}
\dok
(a) $\Leftrightarrow $ (b) is true since the topology on the set 2 is discrete and
the convergence of sequences in Tychonov products is the pointwise convergence
(see \cite{Eng85}).

(b) $\Rightarrow $ (c). By (b), for each $\alpha \in \kappa$ there is $k\in \omega$ such that
$X_n \cap \{ \alpha \} = X\cap \{ \alpha \}$, for each $n\geq k$, which, by Fact \ref{T1244}(b), implies that
$\liminf _{n\in \omega}X_n \cap \{ \alpha \}= \limsup _{n\in \omega}X_n \cap \{ \alpha \}= X\cap \{ \alpha \}$.
This holds for all $\alpha \in \kappa$ so (c) is true.

(c) $\Rightarrow $ (b). Assuming (c), in order to prove (b) we take $\alpha \in \kappa$.
If $\alpha \in X$ then, by (c), there is $k\in \omega$ such that for each $n\geq k$ we have
$\alpha \in X_n$, that is $x_n(\alpha ) = 1 = x(\alpha )$.
If $\alpha \in \kappa \setminus X = \bigcup _{k\in \omega}\bigcap _{n\geq k}\kappa \setminus X_n$ then
there is $k\in \omega$ such that for each $n\geq k$ we have
$\alpha \in \kappa \setminus X_n$, that is $x_n(\alpha ) = 0 = x(\alpha )$ and (b) is true.

The Cantor space $2^\omega$ is sequential, since it is metrizable (see \cite{Eng85}).

Let $\kappa > \omega$ and let $A\subset 2^{\kappa }$ be the family of characteristic functions of
at most countable subsets of $\kappa$. By (a) and since the limit superior
of a sequence of countable sets is countable, the set $A$ is sequentially
closed. But $A$ is dense in $2^{\kappa }$ and, hence, not closed.
Thus $\langle 2^\kappa, \tau_{\mathrm C}\rangle $ is not a sequential space.
\kdok
Let the convergence $\lambda _{\mathrm{s}} $ on the power set $P(\kappa )$ be defined by
$$
\lambda_{\mathrm s}(\langle X_n \rangle)=\left\{
                \begin{array}{cl}
                  \{ X \} & \mbox{ if }{\liminf X_n =\limsup X_n =X ,} \\
                  \emptyset & \mbox{ if } {\liminf X_n <\limsup X_n .}
                \end{array}
              \right.
$$
\begin{fac} \rm \label{T1255}
Let $f : 2^\kappa\rightarrow P(\kappa )$ be the bijection given by $f(x)=x^{-1}[ \{ 1\} ]$. Then

(a) $\tau _{{\mathrm C}}^{P(\kappa )} = \{ f[O] :O\in \tau _{{\mathrm C}} \}$
is a topology on the power set algebra $P(\kappa)$;

(b) $f: \langle 2^\kappa , \tau _{{\mathrm C}} \rangle \rightarrow
\langle P(\kappa ), \tau _{{\mathrm C}}^{P(\kappa )} \rangle$ is a homeomorphism;

(c) $\lambda _{\mathrm{s}}  = \lim_{\tau _{{\mathrm C}}^{P(\kappa )}}= \lim _{{\mathcal O}_{\lambda _{\mathrm{s}} }}$, thus
$\lambda _{\mathrm{s}}  $ is a topological convergence;

(d) ${\mathcal O}_{\lambda _{\mathrm{s}} }= \tau _{{\mathrm C}}^{P(\kappa )}$ iff $\kappa = \omega$. If $\kappa > \omega$, then
$\tau _{{\mathrm C}}^{P(\kappa )} \varsubsetneq {\mathcal O}_{\lambda _{\mathrm{s}} }$.
\end{fac}
\dok
(a) and (b) are evident. Let us prove (c).
By Fact \ref{T1245}, $X\in \lambda _{\mathrm{s}}(\langle X_n \rangle)$ iff
$\langle x_n \rangle \rightarrow _{\tau _{\mathrm{C}}} x$ which is, by (b), equivalent to
$X\in \lim_{\tau _{{\mathrm C}}^{P(\kappa )}}(\langle X_n \rangle)$.
Now, the second equality follows from Fact \ref{T1201}(h).
(d) follows from Fact \ref{T1245} and Fact \ref{T1201}(c).
\kdok

\noindent
The convergence $\lambda_{\mathrm s}$ on the power set algebras is generalized for an arbitrary complete Boolean algebra
${\mathbb B}$ defining the {\bf algebraic convergence} $\lambda_{\mathrm s}$ on ${\mathbb B}$ by
$$
\lambda_{\mathrm s}(\langle x_n \rangle)=\left\{
                \begin{array}{cl}
                  \{ x \} & \mbox{ if }{\liminf x_n =\limsup x_n =x ,} \\
                  0       & \mbox{ if } {\liminf x_n <\limsup x_n ,}
                \end{array}
              \right.
$$
where
$\liminf x_n  = \textstyle \bigvee _{k\in \omega}\bigwedge _{n\geq k}x_n$ and
$\limsup x_n  = \textstyle\bigwedge _{k\in \omega}\bigvee _{n\geq k}x_n $.
By Fact \ref{T1201}, there is the maximal topology ${\mathcal O}_{\lambda_{\mathrm s}}$
on  ${\mathbb B}$ such that $\lambda _{\mathrm{s}}  \leq \lim _{{\mathcal O}_{\lambda_{\mathrm s}}}$,
called the {\bf sequential topology},  traditionally denoted by $\tau _{\mathrm s}$.
It played a significant role in the solution of
von Neumann's \cite{Sco81} and Maharam's Problem  \cite{Mah47} solved by Talagrand \cite{Tal1,Tal2}
(see also papers of  Balcar, Gl\'{o}wczy\'{n}ski and Jech \cite{BGJ98};
Balcar, Jech and Paz\'ak  \cite{BJP05}; Balcar and Jech \cite{BJ06}; Farah \cite{Far04}; Todor\v cevi\'c \cite{Tod04}
and Veli\v ckovi\'c  \cite{Vel05}).

It is known that the convergence $\lambda_{\mathrm s}$ is weakly-topological. Namely we have

\begin{fac} \rm \label{T1246}
Let ${\mathbb B}$ be a complete Boolean algebra. Then

(a) $\lambda_{\mathrm s}$ is a Hausdorff convergence satisfying (L1) and (L2);

(b) $\lambda_{\mathrm s}$ is a weakly-topological convergence.
\end{fac}
\dok
Clearly, $\lambda_{\mathrm s}$ is a Hausdorff convergence and
satisfies (L1). Since for each $x,y\in {\mathbb B}$,  $y \prec x$
implies $\liminf x \leq \liminf y \leq  \limsup y \leq \limsup x$,
it satisfies (L2). (b) follows from (a) and Fact \ref{T1220}.
\kdok
\noindent
By Fact \ref{T1255}(c), on each power set algebra the convergence $\lambda_{\mathrm s}$ is topological.
In general, by Fact \ref{T1201}(h) and Theorem 2 of \cite{KuPa07} we have
\begin{te}\rm\label{T1256}
For each c.B.a.\ ${\mathbb B}$ the following conditions are equivalent:

(a) $\lambda_{\mathrm s}$ is a topological convergence;

(b) $\lambda _{\mathrm{s}}  =  \lim _{{\mathcal O}_{\lambda _{\mathrm{s}} }}$;

(c) The algebra ${\mathbb B}$ is $(\omega , 2)$-distributive;

(d) Forcing by ${\mathbb B}$ does not produce new reals.
\end{te}
\noindent
If an algebra ${\mathbb B}$ is not $(\omega , 2)$-distributive
but
$$
\forall x \in {\mathbb B}^\omega ~ \exists y \prec x~ \forall z \prec y ~ \limsup z =\limsup y ,\eqno{( \hbar )}
$$
then the convergence $\lim_{{\mathcal O}_{\lambda_{\mathrm s}}}$ is characterized
in the following way (see \cite{KuPa07}).
\begin{te} \rm \label{T1257}
If a complete Boolean algebra ${\mathbb B}$ satisfies condition $(\hbar)$,
then for each sequence $x\in {\mathbb B}^\omega $ and $a\in {\mathbb B}$ we have:
$a\in \lim_{{\mathcal O}_{\lambda_{\mathrm s}}}(x)
\Leftrightarrow a_x =b_x=a$, where
$$
a_x  =  \textstyle \bigwedge_{A\in [\omega]^\omega}\bigvee_{B\in [A]^\omega} \bigwedge_{n\in B} x_n \;\;\mbox{ and }\;\;
b_x  =  \textstyle \bigvee_{A\in [\omega]^\omega} \bigwedge_{B\in [A]^\omega}\bigvee_{n\in B} x_n .
$$
The implication ``$\Rightarrow$" holds in each c.B.a.
\end{te}
We note that, by \cite{KuPa07}, condition ($\hbar$) is related to the cellularity of complete Boolean algebras:
$\mathfrak{t}\mbox{-cc}\Rightarrow (\hbar ) \Rightarrow \mathfrak{s}\mbox{-cc}$. By
\cite{KuTo}, $\{ \kappa \in {\mathrm{Card}} :  \kappa \mbox{-cc }
\Rightarrow (\hbar ) \}$ is either  $[0, {\mathfrak h})$  or $ [0,
{\mathfrak h}]$ and $\{ \kappa \in {\mathrm{Card}} : (\hbar
)\Rightarrow \kappa \mbox{-cc }    \} =[{\mathfrak s}, \infty )$.

\section{The Aleksandrov cube and the convergences $\lambda _{\mathrm{ls} }$ and
                              $\lambda _{\mathrm{li} }$}

We remind the reader that the Aleksandrov cube  of weight $\kappa$,
here denoted by $\langle 2^\kappa , \tau _A \rangle$, is the Tychonov product of $\kappa$-many
copies of the two-point space $2=\{ 0, 1 \}$ with the topology
${\mathcal O}_A= \{ \emptyset ,  \{ 0\} , \{ 0,1 \}\}$. It is an universal T$_0$ space of weight $\kappa$
(see \cite{Eng85}).

\begin{fac} \rm \label{T1254}
(a) Let $\langle x_n \rangle$ be a sequence in $2^\kappa$ and $x\in 2^\kappa$.
Then $\langle x_n \rangle \rightarrow_{\tau _A } x$ iff
\begin{equation}\label{EQ1243}
  \textstyle  \limsup _{n\in \omega}X_n \subset X,
\end{equation}
where $X_n = x_n ^{-1}[\{ 1\}]$, for $n\in \omega $, and  $X = x ^{-1}[\{ 1\}]$.

(b) $\langle 2^\kappa, \tau_A \rangle $ is a sequential space iff $\kappa=\omega$.
\end{fac}
\dok
(a) In the space $\langle 2,{\mathcal O}_A \rangle$ the point 0
is isolated and the only neighborhood of the point 1 is $\{ 0,1 \}$
so, a sequence $\langle a_n : n\in \omega  \rangle$ converges to a
point $a$ iff $a=1$, or $a=0$ and there is $k\in \omega $ such that
$a_n=0$, for all $n\geq k$. Now as in Section 3 we conclude that, in
the space $\langle 2^\kappa , \tau _A \rangle$, $\langle x_n \rangle
\rightarrow _{\tau _A} x$ iff for each $\alpha <\kappa$, $\langle
x_n (\alpha ) \rangle \rightarrow _{{\mathcal O} _A} x(\alpha )$ iff
$$
\forall \alpha <\kappa \;\;
\Big[ x(\alpha )=1 \; \lor \;
\Big( x(\alpha )=0 \; \land \; \exists k\in \omega\;\; \forall n\geq k \;\; x_n(\alpha ) =0\Big)\Big]
$$
iff for each $\alpha < \kappa $ we have $\alpha \in X \lor \neg \forall k\in \omega \; \exists n\geq k \; \alpha \in X_n$,
that is $\alpha \in \limsup X_n \Rightarrow \alpha \in X$.

(b)
($\Leftarrow$) $\langle 2^\omega , \tau _A \rangle$ is a first countable and, consequently, a sequential space.

($\Rightarrow$) Let $\kappa > \omega$. The set $S=\{ x\in 2^\kappa : |x^{-1}[\{ 0 \}]| \leq \aleph _0 \}$
is dense in the space $\langle 2^\kappa , \tau _A \rangle$ and, hence, it is not closed. In order to show that $S$ is
sequentially closed we take a sequence $\langle x_n : n\in \omega  \rangle$ in $S$ and show that
$\lim _{\tau _A}(\langle x_n \rangle) \subset S$. The corresponding sets
$X_n = x_n ^{-1}[\{ 1\}]$, $n\in \omega $, are co-countable subsets of $\kappa$, thus $X_n = \kappa \setminus C_n$,
where $C_n \in [\kappa ]^{\leq \aleph _0}$ and the set
$\limsup X_n=\kappa \setminus \bigcup _{k\in \omega }\bigcap _{n\geq k}C_n$ is co-countable as well.
By (a), if $x\in \lim _{\tau _A}(\langle x_n \rangle)$, then $\limsup X_n \subset X$, which means that $X$
is a co-countable set and, consequently, $x\in S$.
\kdok
Let the
convergence $\lambda _{\mathrm{ls} }$ on $P(\kappa )$ be defined by
$$
\lambda_{\mathrm{ls} }(\langle X_n \rangle)=(\limsup X_n )\upar.
$$
\begin{te} \rm \label{T1260}
Let $f : 2^\kappa\rightarrow P(\kappa )$ be the bijection given by $f(x)=x^{-1}[ \{ 1\} ]$. Then

(a) $\tau _{{\mathrm A}}^{P(\kappa )} = \{ f[O] :O\in \tau _{{\mathrm A}} \}$ is a
topology on $P(\kappa)$;

(b) $f: \langle 2^\kappa , \tau _{{\mathrm A}} \rangle
\rightarrow \langle P(\kappa ), \tau _{{\mathrm A}}^{P(\kappa )} \rangle$
is a homeomorphism;

(c) $\lambda _{\mathrm{ls} } = \lim_{\tau _{{\mathrm A}}^{P(\kappa )}}= \lim _{{\mathcal O}_{\lambda _{\mathrm{ls} }}}$
    and $\lambda _{\mathrm{ls} } $ is a topological convergence.

(d) ${\mathcal O}_{\lambda _{\mathrm{ls} }}= \tau _{{\mathrm A}}^{P(\kappa )}$ iff $\kappa = \omega$.
If $\kappa > \omega$, then $\tau _{{\mathrm A}}^{P(\kappa )} \varsubsetneq {\mathcal O}_{\lambda _{\mathrm{ls} }}$.

(e) ${\mathcal O}_{\lambda _{\mathrm{ls} }}\not\subset \tau _{{\mathrm C}}^{P(\kappa )}$, if $\kappa > \omega$.
\end{te}
\dok
(a) and (b) are evident.
(c) and (d) follow from Fact \ref{T1254} and Fact \ref{T1201}(h).

(e) As in Fact \ref{T1254} we consider the set $F= \{ \kappa
\setminus C : C \in [\kappa ]^{\leq \aleph _0 }\}$, which is dense
in the space $\langle P(\kappa ),\tau _{{\mathrm C}}^{P(\kappa )}
\rangle$ and, hence $P(\kappa )\setminus F \not\in \tau _{{\mathrm
C}}^{P(\kappa )}$. If $\langle X_n \rangle $ is a sequence in $F$,
where $X_n = \kappa \setminus C_n$, then $\limsup X_n=\kappa
\setminus \bigcup _{k\in \omega }\bigcap _{n\geq k}C_n \in F$ and,
clearly, $\lambda _{\mathrm{ls} }(\langle X_n \rangle )= (\limsup
X_n)\upar \subset F$, thus $u_{\lambda _{\mathrm{ls} }}(F)=F$. By
(c), $\lambda _{\mathrm{ls} }$ satisfies (L1) and (L2) so, by Fact
\ref{T1201}(d), $P(\kappa )\setminus F \in {\mathcal O}_{\lambda _{\mathrm{ls} }}$.
\kdok

Now we generalize this for an arbitrary complete Boolean algebra ${\mathbb B}$ defining
the convergence $\lambda_{\mathrm{ls} }$ by
$$
\lambda_{\mathrm{ls} }(\langle x_n \rangle)=(\limsup x_n)\upar
$$
and Fact \ref{T1201} provides the topology ${\mathcal O}_{\lambda_{\mathrm{ls} }}$ on  ${\mathbb B}$.
We will also consider the dual convergence $\lambda_{\mathrm{li} }$
on ${\mathbb B}$ defined by
$
\lambda_{\mathrm{li} }(\langle x_n \rangle)=(\liminf x_n)\downar
$
and the corresponding topology ${\mathcal O}_{\lambda_{\mathrm{li} }}$.

If $\lambda _1$ and $\lambda _2 $ are convergences, by $\lambda _1 \cap \lambda _2$ we will denote
the convergence defined by $(\lambda _1 \cap \lambda _2) (x) =\lambda _1 (x) \cap \lambda _2 (x)$.
Similarly to Fact \ref{T1246} we have

\begin{te} \rm \label{T1261}
Let ${\mathbb B}$ be a complete Boolean algebra. Then

(a) $\lambda _{\mathrm{ls} }$ and $\lambda _{\mathrm{li} }$ are non-Hausdorff convergences satisfying (L1) and (L2);

(b) $\lambda_{\mathrm s} = \lambda_{\mathrm{ls} } \cap \lambda_{\mathrm{li} }$ and,
consequently, $\lambda_{\mathrm s}\leq \lambda_{\mathrm{ls}} , \lambda_{\mathrm{li}}$;

(c) ${\mathcal O}_{\lambda_{\mathrm{ls} }}, {\mathcal O}_{\lambda_{\mathrm{li} }}
     \subset {\mathcal O}_{\lambda_{\mathrm s}}$;

(d) $\lambda _{\mathrm{ls} }^* \leq \lim _{{\mathcal O}_{\lambda _{\mathrm{ls} }}}$ and
    $\lambda _{\mathrm{li} }^* \leq \lim _{{\mathcal O}_{\lambda _{\mathrm{li} }}}$.

(e) $\lambda_{\mathrm s}^* = \lambda_{\mathrm{ls} }^* \cap \lambda_{\mathrm{li} }^*$ and,
consequently, $\lambda_{\mathrm s}^*\leq \lambda_{\mathrm{ls}}^* , \lambda_{\mathrm{li}}^*$.
\end{te}

\dok
(a) Since $a\in a\upar =(\limsup \langle a \rangle )\upar =
\lambda _{\mathrm{ls} }(\langle a \rangle ) $, for each $a\in
{\mathbb B}$, $\lambda _{\mathrm{ls} }$ satisfies (L1) and it is not
Hausdorff because $0,1\in \lambda _{\mathrm{ls} }(\langle 0 \rangle
) $. For a proof of (L2) note that $y\prec x$ implies $\limsup y
\leq \limsup x$ so we have $(\limsup x) \upar \subset (\limsup y)
\upar $, that is $\lambda _{\mathrm{ls} }(x) \subset \lambda
_{\mathrm{ls} }(y)$.

(b) If $a\in \lambda_{\mathrm s}(x)$, then $a=\limsup x \in (\limsup
x)\upar = \lambda _{\mathrm{ls} }(x)$ and, similarly, $a\in \lambda
_{\mathrm{li} }(x)$. Conversely, if $a\in \lambda_{\mathrm{ls} }(x)
\cap \lambda_{\mathrm{li} }(x)$, then $\limsup x \leq a \leq \liminf
x$, which implies $\limsup x = \liminf x =a$, that is $a\in
\lambda_{\mathrm s}(x)$.

(c) follows from (b) and Fact \ref{T1201}(g). (d) follows from Fact
\ref{T1247}.

(e) By (b) and by the minimality of $\lambda ^*$ (see Fact \ref{T1247})
we have $\lambda _{\mathrm{s}}^* \leq \lambda
_{\mathrm{ls}}^* ,\lambda _{\mathrm{li}}^*$. So, it remains to be
proved that $\lambda _{\mathrm{ls}}^* \cap \lambda _{\mathrm{li}}^*
\leq\lambda _{\mathrm{s}}^*$. Let $x\in {\mathbb B}^\omega $ and
$a\in \lambda _{\mathrm{ls}}^*(x) \cap \lambda _{\mathrm{li}}^*(x)$.
If $y\prec x$, then there is $z\prec y$ such that $a\geq \limsup z$
and there is $t\prec z $ such that $a\leq \liminf t$. But then
$\limsup t \leq \limsup z \leq a \leq \liminf t$, which implies $a
\in \lambda _{{\mathrm s}}(x)$. Thus for each $y\prec x$ there is
$t\prec y$ such that $a \in \lambda _{{\mathrm s}}(x)$, that is $a \in \lambda _{{\mathrm s}}^*(x)$.
\kdok
By the previous theorem and Fact \ref{T1246}, the relations between the convergences considered in this paper are presented in
 the following diagram.
\vspace{5mm}
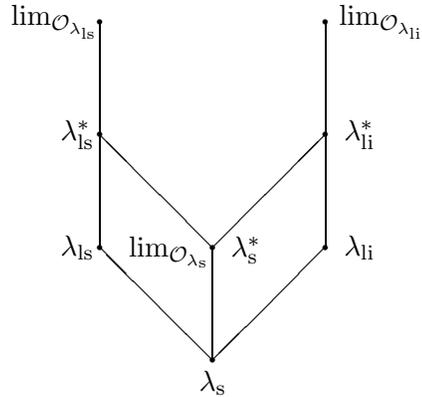
\begin{center}
\unitlength .3mm 
\linethickness{0.4pt}
\ifx\plotpoint\undefined\newsavebox{\plotpoint}\fi 
\begin{picture}(200,176)(0,0)
\put(150,175){\line(0,-1){50}}
\put(150,125){\line(0,-1){50}}
\put(150,75){\line(-1,-1){50}}
\put(100,25){\line(-1,1){50}}
\put(50,75){\line(0,1){100}}
\put(100,25){\line(0,1){50}}
\put(100,75){\line(1,1){50}}
\put(100,75){\line(-1,1){50}}
\put(150,175){\circle*{2}}
\put(150,125){\circle*{2}}
\put(150,75){\circle*{2}}
\put(100,25){\circle*{2}}
\put(100,75){\circle*{2}}
\put(50,125){\circle*{2}}
\put(50,175){\circle*{2}}
\put(50,75){\circle*{2}}
\put(175,175){\makebox(0,0)[cc]{$\lim_{{\mathcal O}_{\lambda_{{\mathrm{li}}}}}$}}
\put(165,125){\makebox(0,0)[cc]{$\lambda^*_{{\mathrm{li}}}$}}
\put(40,125){\makebox(0,0)[cc]{${\lambda^*_{{\mathrm{ls}}}}$}}
\put(100,15){\makebox(0,0)[cc]{$\lambda_{{\mathrm{s}}}$}}
\put(40,75){\makebox(0,0)[cc]{$\lambda_{{\mathrm{ls}}}$}}
\put(165,75){\makebox(0,0)[cc]{$\lambda_{{\mathrm{li}}}$}}
\put(92,73){\makebox(0,0)[cc]{$\lim_{{\mathcal O}_{\lambda_{{\mathrm{s}}}}}\;\;\lambda^*_{{\mathrm{s}}}$}}
\put(30,175){\makebox(0,0)[cc]{$\lim_{{\mathcal O}_{\lambda_{{\mathrm{ls}}}}}$}}
\end{picture}
\end{center}
In the sequel we will use the following characterization, where
the families of closed sets corresponding to the topologies
${\mathcal O}_{\lambda_{\mathrm{ls}}}$ and ${\mathcal O}_{\lambda_{\mathrm{li}}}$ are denoted by
${\mathcal F}_{\lambda_{\mathrm{ls}}}$ and ${\mathcal F}_{\lambda_{\mathrm{li}}}$ respectively.

\begin{te}  \rm \label{T4105}
Let ${\mathbb B}$  be a complete Boolean algebra. Then

\noindent
(I) For a set $F \subset {\mathbb B}$ the following conditions are equivalent:

(a) $F\in {\mathcal F}_{\lambda_{\mathrm{ls}}}$;

(b)  $F$  is upward closed and $\limsup x\in F$, for each sequence $x\in F^\omega$;

(c) $F$ is upward closed and  $\bigwedge_{n \in \omega} x_n \in F$, for each decreasing $x\in F^\omega$.

\noindent
(II) For a set $F \subset {\mathbb B}$ the following conditions are equivalent:

(a) $F\in {\mathcal F}_{\lambda_{\mathrm{li}}}$;

(b)  $F$  is downward closed and $\liminf x\in F$, for each sequence $x\in F^\omega$;

(c) $F$ is downward closed and  $\bigvee_{n \in \omega} x_n \in F$,
for each increasing $x\in F^\omega$.

\noindent
(III) The mapping $h:\langle {\mathbb B},{\mathcal O_{\lambda_{\mathrm{ls}}}}\rangle \rightarrow
                     \langle {\mathbb B},{\mathcal O_{\lambda_{\mathrm{li}}}}\rangle$
                     given by $h(b)=b'$, for each $b\in  {\mathbb B}$, is a homeomorphism.
\end{te}
\dok
We prove (I). The proof of (II) is dual.

(a) $\Rightarrow$ (b). Let $X\setminus F \in {\mathcal
O}_{\lambda_{\mathrm{ls}}}$. Then, by Theorems \ref{T1261} and Fact
\ref{T1201}(d) we have $F=u_{\lambda_{\mathrm{ls}}}(F)=\bigcup
_{x\in F^\omega }(\limsup x )\upar$ and, hence, $F$ is upward
closed. Also, if $x\in F^\omega$, then $\limsup x \in (\limsup
x)\upar \subset F$.

(b) $\Rightarrow$ (c). If $x\in F^\omega $ is a decreasing sequence, $ \bigwedge_{n \in \omega} x_n = \limsup x \in F$.

(c) $\Rightarrow$ (a). Assuming (c), by Fact \ref{T1201}(d) we show
that $u_{\lambda_{\mathrm{ls}}}(F)=F$. If $b \in
u_{\lambda_{\mathrm{ls}}}(F)$, then there is $x \in F^\omega $ such
that $b \geq \limsup x$. Since the set $F$ is upward closed and
$x_n\in F$, for $k\in \omega$ we have $y_k=b \vee \bigvee_{n\geq k}
x_n\in F$ and, clearly, $y=\langle y_k : k\in \omega  \rangle$ is a
decreasing sequence. So, $F\ni \bigwedge_{k\in \omega}y_k
=\bigwedge_{k\in \omega}(b\vee \bigvee_{n\geq k}x_n) =b \vee
\bigwedge_{k\in \omega} \bigvee_{n\geq k} x_n = b \vee \limsup x
=b$.

(III) $h$ is a bijection and for a proof of its continuity
we take $F\in {\mathcal F}_{\lambda_{\mathrm{li}}}$ and show that
$h^{-1}[F]=\{ b' : b\in F \} \in {\mathcal
F}_{\lambda_{\mathrm{ls}}}$. If $a\geq b' \in h^{-1}[F]$, then
$a'\leq b\in F$ and, by (II), $a'\in F$, which implies $a\in
h^{-1}[F]$. Thus the set $ h^{-1}[F]$ is upward closed. Let $\langle
x_n \rangle$ be a decreasing sequence in $h^{-1}[F]$. Then $\langle
x_n' \rangle$ is an increasing sequence in $F$ and, by (II) again,
$\bigvee _{n\in \omega} x_n' = (\bigwedge _{n\in \omega} x_n)'\in
F$, which implies $\bigwedge _{n\in \omega} x_n \in h^{-1}[F]$. By
(I), $h^{-1}[F] \in {\mathcal F}_{\lambda_{\mathrm{ls}}}$. The proof
that $h$ is closed is similar.$\Box$

\section{The algebras with $\lambda _{\mathrm{ls}}$ topological}

In this section we prove the following characterization of complete Boolean algebras
on which the convergences  $\lambda _{\mathrm{ls}}$ and $\lambda _{\mathrm{li}}$ are topological.

\begin{te}  \rm \label{T1280}
For each complete Boolean algebra ${\mathbb B}$ the following conditions are equivalent:

(a) $\lambda_{\mathrm{ls}}$ is a topological convergence;

(b) $\lambda_{\mathrm{li}}$ is a topological convergence;

(c) ${\mathbb B}$ is an $(\omega,2)$-distributive algebra;

(d) Forcing by ${\mathbb B}$ does not produce new reals.
\end{te}
The following three lemmas will be used in our proof.
\begin{lem}  \rm \label{T4116}
Let ${\mathbb B}$  be a complete Boolean algebra. Then

(a) For each $a \in {\mathbb B}$ the function
$f_a :\langle {\mathbb B},{\mathcal O_{\lambda_{\mathrm{ls}}}}\rangle \rightarrow
\langle {\mathbb B},{\mathcal O_{\lambda_{\mathrm{ls}}}}\rangle$
defined by $f_a (x)=x \wedge a$ is continuous;

(b) $\Lim_{{\mathcal O}_{\lambda_{\mathrm{ls}}}}\neq{\lambda}_{\mathrm{ls}}$ iff
there is a sequence $x$ in ${\mathbb B}$  such that
$0 \in \Lim_{{\mathcal O}_{\lambda_{\mathrm{ls}}}}(x) \setminus {\lambda}_{\mathrm{ls}}(x)$.

(c) If $x,y\in {\mathbb B}^\omega$ and $x_n \leq y_n$, for each $n \in \omega$,
then $\Lim_{{\mathcal O}_{\lambda_{\mathrm{ls}}}} (y) \subset \Lim_{{\mathcal O}_{\lambda_{\mathrm{ls}}}}(x)$.
\end{lem}
\dok
(a) By Theorem \ref{T4105} we show that for a closed set
$F\subset {\mathbb B}$ the set $f_a^{-1}[F]=\{x \in {\mathbb B}: x
\wedge a \in F\}$ is upward closed and contains the infimum of each
decreasing sequence in $f_a^{-1}[F]$. First, if $x_1\geq x\in
f_a^{-1}[F]$, then $x_1 \wedge a \geq x \wedge a \in F$ and, since
$F$ is upward closed, $x_1\wedge a \in F$, that is $x_1 \in
f_a^{-1}[F]$. Second, if $\langle x_n \rangle $ is a decreasing
sequence in $f_a^{-1}[F]$, then $\langle x_n\wedge a \rangle $ is a
decreasing sequence in $F$ and, since $F$ is closed,
$\bigwedge_{n\in \omega} x_n\wedge a \in F$, thus $\bigwedge_{n\in
\omega} x_n  \in f_a^{-1}[F]$.

(b) Let $y\in {\mathbb B}^\omega$ and $b \in \Lim_{{\mathcal O}_{\lambda_{\mathrm{ls}}}}(y) \setminus {\lambda}_{\mathrm{ls}}(y)$. Then
$ \limsup y \not \leq b$ and, hence, $c=\limsup y \wedge b'>0$.
Let $x= \langle y_n  \wedge c :n \in \omega\rangle$.
Since $c\leq \limsup y $ we have $c =\bigwedge_{k \in \omega} \bigvee_{n\geq k}y_n \wedge c = \limsup x$,
which implies $0 \not \in \lambda_{\mathrm{ls}}(x)$.
Since $b \in \Lim_{{\mathcal O}_{\lambda_{\mathrm{ls}}}}(y)$ and, by (a), the function
$f_c:{\mathbb B}\rightarrow {\mathbb B}$ defined by  $f_c(t)=t\wedge c$ is continuous,
we have $0=b\wedge c  = f_c (b)\in  \Lim_{{\mathcal O}_{\lambda_{\mathrm{ls}}}}(\langle f_c(y_n) \rangle)
= \Lim_{{\mathcal O}_{\lambda_{\mathrm{ls}}}}(x)$.

(c) Let $a\in \Lim_{{\mathcal O}_{\lambda_{\mathrm{ls}}}} (y)$ and $a\in O \in {\mathcal O}_{\lambda_{\mathrm{ls}}}$.
Then there is $n_0\in \omega $ such that for each $n\geq n_0$ we have $y_n\in O$, thus, since by Theorem \ref{T4105}
the set $O$ is downward closed, $x_n \in O$, for $n\geq n_0$. So $a\in \Lim_{{\mathcal O}_{\lambda_{\mathrm{ls}}}} (x)$.
\kdok

\noindent
If $x\in {\mathbb B}^\omega $, then $\tau
_x= \{ \langle \check{n}, x_n \rangle : n\in \omega \}$ is the
corresponding ${\mathbb B}$-name for a subset of $\omega$ and, by
Lemmas 2 and 6 of \cite{KuPa07},
\begin{eqnarray*}
\liminf x & = & \| \check{\omega} \subset ^* \tau_x \| ;\\
\limsup x & = & \| |\tau_x |=\check{\omega}\| ;\\
a_x       & = & \|\forall A\in(([\omega]^\omega)^V)^{\check{~}}\, \exists B\in(([A]^\omega)^V)^{\check{~}} \,
                B\subset^*\tau_x \| ;\\
b_x       & = & \|\exists A\in(([\omega]^\omega)^V)^{\check{~}}\, \forall B\in(([A]^\omega)^V)^{\check{~}} \,
                |\tau_x\cap B|=\check{\omega}\| .
\end{eqnarray*}
\begin{lem}\rm\label{T1264}
Let ${\mathbb B}$ be a complete Boolean algebra and $x$ a sequence in ${\mathbb B}$. Then

(a) $\liminf x \leq a_x \leq b_x \leq \limsup x$;

(b) If ${\mathbb B}$ is $(\omega , 2)$-distributive, then
    $a_x= \liminf x$ and $b_x=\limsup x$;

(c) $b_x=\bigvee_{y\prec x}\bigwedge_{z\prec y} \bigvee_{m \in \omega}z_m$.
\end{lem}

\dok
(a) This is Lemma 7 of \cite{KuPa07}.

(b) Let ${\mathbb B}$ be $(\omega,2)$-distributive. By (a), it is
sufficient to show that $\limsup x \leq b_x$, that is $1\Vdash
|\tau_x |=\check{\omega} \Rightarrow \exists
A\in(([\omega]^\omega)^V)^{\check{~}}\, \forall
B\in(([A]^\omega)^V)^{\check{~}} \, |\tau_x\cap B|=\check{\omega}$.
Let $G$ be a ${\mathbb B}$-generic filter over $V$ and let $|(\tau_x
)_G |=\omega$. Then, by the $(\omega,2)$-distributivity we have
$(\tau_x )_G \in ([\omega]^\omega)^V$ and $A=(\tau_x )_G$ is as
required. Thus $b_x=\limsup x$. The proof of $a_x= \liminf x$ is
similar.

(c) Clearly we have
$\bigvee_{y\prec x}\bigwedge_{z\prec y} \bigvee_{m \in \omega} z_m
=\bigvee_{f\in \omega^{\uparrow \omega}}\bigwedge_{z\prec x\circ f} \bigvee_{m \in \omega} z_m=$
$\bigvee_{f\in \omega^{\uparrow \omega}}\bigwedge_{g\in \omega^{\uparrow \omega}} \bigvee_{m \in \omega} x_{f(g(m))}
=\bigvee_{f\in \omega^{\uparrow \omega}}\bigwedge_{g\in \omega^{\uparrow \omega}} \bigvee_{n \in f[g[\omega]]}x_n$
and we prove that in each generic extension $V_{{\mathbb B}}[G]$ conditions
\begin{equation}\label{EQ1250}
\exists A \in ([\omega ]^\omega )^V \; \forall B \in  ([A]^\omega )^V \;
B\cap (\tau _x)_G\neq \emptyset \mbox{ and }
\end{equation}
\begin{equation}\label{EQ1251}
\exists f\in (\omega^{\uparrow \omega})^V \; \forall g\in (\omega^{\uparrow \omega})^V \;
f[g[\omega ]]\cap (\tau _x)_G\neq \emptyset
\end{equation}
are equivalent.
Let (\ref{EQ1250}) hold and let $f_A$ be the increasing enumeration of the set $A$. Then
$f_A \in (\omega^{\uparrow \omega})^V$ and for any
$g\in (\omega^{\uparrow \omega})^V$ we have $f_A[g[\omega ]]\in ([A]^\omega )^V$ thus, by the assumption,
$f_A [g[\omega ]]\cap (\tau _x)_G\neq \emptyset$.

Let (\ref{EQ1251}) hold. Then $A=f[\omega ]\in ([\omega ]^\omega )^V$ and, if $B \in  ([A]^\omega )^V$, then
$f^{-1}[B]\in ([\omega ]^\omega )^V $ and $g_{f^{-1}[B]}\in (\omega^{\uparrow \omega})^V$, where $g_{f^{-1}[B]}$
is the increasing enumeration of the set $f^{-1}[B]$. By the assumption we have
$f[g_{f^{-1}[B]}[\omega ]]\cap (\tau _x)_G\neq \emptyset$ and, since $f[g_{f^{-1}[B]}[\omega ]]=f[f^{-1}[B]]=B$
(because $B\subset f[\omega ]$), we have $B\cap (\tau _x)_G\neq \emptyset$.
\kdok

A sequence $x$ in a c.B.a.\ ${\mathbb B}$ will be called {\bf lim\,sup-stable} ({\bf lim\,inf-stable}, respectively)
iff $\limsup y = \limsup x$  ($\liminf y = \liminf x$ respectively), for each subsequence $y$ of $x$.

\begin{lem}\rm \label{T1233a}
Let $x=\langle x_n: n \in \omega\rangle $ be a sequence in a c.B.a.\ ${\mathbb B}$.

(a) If $x$ is a lim\,sup-stable sequence, then in the space
$\langle {\mathbb B},{\mathcal O}_{\lambda_{\mathrm{ls}}} \rangle$ we have
\begin{equation}\label{EQ1233c} \textstyle
\overline{\{x_n:n \in \omega\}}=(\limsup x)\upar \cup \;
\bigcup_{n\in \omega} x_n\upar  ;
\end{equation}

(b) If $x$ is a lim\,inf-stable sequence, then in the space
$\langle {\mathbb B},{\mathcal O}_{\lambda_{\mathrm{li}}} \rangle$ we have
\begin{equation}\label{EQ1233f} \textstyle
\overline{\{x_n:n \in \omega\}}=(\liminf x)\downar \cup \;
\bigcup_{n\in \omega} x_n\downar  .
\end{equation}
\end{lem}

\dok
We prove (a) and the proof of (b) is dual.
Let $X=\{ x_n : n\in \omega \}$. First we prove that
\begin{equation}\label{EQ1233a} \textstyle
u_{\lambda _{{\mathrm{ls}}}}(X)=(\limsup x)\upar \cup \;
\bigcup_{n\in \omega} x_n\upar  .
\end{equation}
Since $(\limsup x)\upar = \lambda _{{\mathrm{ls}}}(\langle x_n: n
\in \omega\rangle )$ and $x_n \upar = \lambda
_{{\mathrm{ls}}}(\langle x_n , x_n , \dots \rangle )$, for each
$n\in \omega $, the inclusion ``$\supset $" in (\ref{EQ1233a}) is
proved. By Theorems \ref{T1261}(a) and Fact \ref{T1206} we have
$\overline{X} = u_{\lambda _{{\mathrm{ls}}}}^{\omega _1}(X)\supset
u_{\lambda _{{\mathrm{ls}}}}(X)$ and the inclusion ``$\supset $" in
(\ref{EQ1233c}) is true as well.

In order to prove the inclusion ``$\subset$" in (\ref{EQ1233a}) we
take $y\in X^\omega$. If $y$ has a constant subsequence, say
$\langle x_n, x_n, \dots \rangle $, then $x_n \leq \limsup y$ and,
hence, $\lambda _{{\mathrm{ls}}}(y)= (\limsup y)\upar \subset x_n
\upar$ and we are done. Otherwise, by Ramsey's Theorem, there is
$H\in [\omega ]^{\omega }$ such that $y\upharpoonright H$ is an
injection. Let the function $f:H\rightarrow \omega$ be defined by
$f(k)=\min\{n\in \omega :y_{k}=x_n\}$. Then for different $k_1, k_2
\in H$ we have $x_{f(k_1)}=y_{k_1}\neq y_{k_2}=x_{f(k_2)}$ and,
hence, $f(k_1)\neq f(k_2)$. Thus $f$ is an injection so, by Ramsey's
Theorem again and since $\omega $ is a well ordering, there is $H_1
\in [H]^{\omega }$ such that $f\upharpoonright H_1$ is an increasing
function. Now we have $y \succ \langle y_k : k\in H_1 \rangle =
\langle x_{f(k)} : k\in H_1\rangle \prec x$ and, since $x$ is a
lim\,sup-stable sequence, $\limsup y \geq \limsup  \langle y_k :
k\in H_1 \rangle = \limsup x$, which implies $\lambda
_{{\mathrm{ls}}}(y)= (\limsup y)\upar \subset (\limsup x)\upar$ and
(\ref{EQ1233a}) is proved.

Now, we  prove that
\begin{equation}\label{EQ1233b}
u_{\lambda _{{\mathrm{ls}}}}(X)=u_{\lambda _{{\mathrm{ls}}}}(u_{\lambda _{{\mathrm{ls}}}}(X)).
\end{equation}
The inclusion ``$\subset $" holds, since $\lambda _{{\mathrm{ls}}}$ satisfies (L1).
In order to prove ``$\supset $" for $y\in u_{\lambda _{{\mathrm{ls}}}}(X)^\omega $
we show that $\lambda _{{\mathrm{ls}}}(y) \subset u_{\lambda _{{\mathrm{ls}}}}(X)$.
By (\ref{EQ1233a}) we have
$$
\forall k\in \omega \;\; (y_k \geq \limsup x \lor \exists n \in \omega \;\; y_k \geq x_n ).
$$
If there exists $G\in [\omega]^\omega$ such that $y_k \geq \limsup
x$, for each $k \in G$, then $\limsup x \leq \limsup \langle y_k :k
\in G\rangle \leq \limsup y$, which implies $\lambda
_{{\mathrm{ls}}}(y)= (\limsup y)\upar \subset (\limsup x)\upar$
$\subset u_{\lambda _{{\mathrm{ls}}}}(X)$.

Otherwise, there is $k_0 \in \omega$ such that for all $k \geq k_0$ there is $n\in \omega$
such that $y_k \geq x_n$. Let $f:\omega \setminus k_0  \rightarrow \omega$ be defined by
$f(k)=\min\{n\in \omega :x_n \leq y_k \}$. Then $x_{f(k)}\leq y_k$, for $k\in \omega \setminus k_0$.

If there are $H_0\in [\omega \setminus k_0]^\omega$ and $n\in
\omega$ such that $f(k)=n$, for each $k\in H_0$, then $\limsup y
\geq \limsup \langle y_k :k \in H_0\rangle \geq x_n$, which implies
$\lambda _{{\mathrm{ls}}}(y)= (\limsup y)\upar \subset x_n \upar
\subset u_{\lambda _{{\mathrm{ls}}}}(X)$. Otherwise, by Ramsey's
Theorem, there is $H_1\in [\omega \setminus k_0]^\omega$ such that
$f\upharpoonright H_1$ is an injection and, by Ramsey's Theorem
again, there exists $H_2 \in [H_1]^\omega$ such that
$f\upharpoonright H_2$  is an increasing mapping. Now $\langle y_k
:k \in H_2\rangle \prec y$, which implies
\begin{equation}\label{EQ1233d}
\limsup\langle y_k:k \in H_2\rangle \leq \limsup y
\end{equation}
and $\langle x_{f(k)}:k \in H_2\rangle \prec x$, which, since $x$ is a lim\,sup-stable sequence, implies
\begin{equation}\label{EQ1233e}
\limsup\langle x_{f(k)}:k \in H_2\rangle = \limsup x .
\end{equation}
Since $x_{f(k)}\leq y_k$ we have $\limsup\langle x_{f(k)}:k \in
H_2\rangle \leq \limsup\langle y_k:k \in H_2\rangle$ and, by
(\ref{EQ1233d}) and (\ref{EQ1233e}), $\limsup x \leq \limsup y$ so
$\lambda _{{\mathrm{ls}}}(y) \subset u_{\lambda
_{{\mathrm{ls}}}}(X)$ again.

Since the convergence $\lambda _{{\mathrm{ls}}}$ satisfies (L1) and
(L2), by Fact \ref{T1206} we have $\overline{X} = u_{\lambda
_{{\mathrm{ls}}}}^{\omega _1}(X)$ and (\ref{EQ1233c}) follows from
(\ref{EQ1233a}) and (\ref{EQ1233b}).
\kdok

\noindent
{\bf Proof of Theorem \ref{T1280}}
(c) $\Leftrightarrow$ (d) is a well known fact (see \cite{Jec97}).

(a) $\Leftrightarrow$ (b). Assuming that $\lambda_{\mathrm{ls}}= \Lim_{{\mathcal O}_{\lambda_{\mathrm{ls}}}}$ we
prove that $\lambda_{\mathrm{li}}= \Lim_{{\mathcal O}_{\lambda_{\mathrm{li}}}}$, that is
$ \Lim_{{\mathcal O}_{\lambda_{\mathrm{li}}}}(x) \subset\lambda_{\mathrm{li}}(x)$, for each sequence
$x$ in ${\mathbb B}$. So, if $a\in \Lim_{{\mathcal O}_{\lambda_{\mathrm{li}}}}(x)$, then, by Theorem \ref{T4105}(III),
we have $a'\in \Lim_{{\mathcal O}_{\lambda_{\mathrm{ls}}}}(\langle x_n' \rangle )
= \lambda_{\mathrm{ls}}(\langle x_n' \rangle )$, that is $a' \geq \limsup x_n'$, which implies
$a\leq \liminf x_n$ and, hence, $a\in \lambda_{\mathrm{li}}(x)$. The proof of the converse is similar.

(a) $\Rightarrow$ (c). If $\lambda_{\mathrm{ls}}$ is a topological
convergence, then $\lambda_{\mathrm{li}}$ is topological as well. By
Theorem \ref{T1261}(c) we have ${\mathcal
O}_{\lambda_{\mathrm{ls}}}, {\mathcal
O}_{\lambda_{\mathrm{li}}}\subset {\mathcal
O}_{\lambda_{\mathrm{s}}}$, and, by Fact \ref{T1200},
$\Lim_{{\mathcal O}_{\lambda_{\mathrm{s}}}}\leq \Lim_{{\mathcal
O}_{\lambda_{\mathrm{ls}}}}, \Lim_{{\mathcal
O}_{\lambda_{\mathrm{li}}}}$ so, since $\lambda_{\mathrm{ls}}$ and
$\lambda_{\mathrm{li}}$ are topological, $\Lim_{{\mathcal
O}_{\lambda_{\mathrm{s}}}}\leq \lambda_{\mathrm{ls}},
\lambda_{\mathrm{li}}$, which, by Theorem \ref{T1261}(b) implies
$\Lim_{{\mathcal O}_{\lambda_{\mathrm{s}}}}\leq
\lambda_{\mathrm{ls}}
\cap\lambda_{\mathrm{li}}=\lambda_{\mathrm{s}}\leq \Lim_{{\mathcal
O}_{\lambda_{\mathrm{s}}}}$. So,
$\lambda_{\mathrm{s}}=\Lim_{{\mathcal O}_{\lambda_{\mathrm{s}}}}$,
that is $\lambda_{\mathrm{s}}$ is a topological convergence and, by
Theorem \ref{T1256}, the algebra ${\mathbb B}$ is
$(\omega,2)$-distributive.

(c) $\Rightarrow$ (a). Suppose that the algebra ${\mathbb B}$ is
$(\omega,2)$-distributive and that $\lambda_{\mathrm{ls}}$ is not a
topological convergence. Then, by Lemma \ref{T4116}(b), there exists
a sequence $x$ in ${\mathbb B}$ such that $0\in \Lim_{{\mathcal
O}_{\lambda_{\mathrm{ls}}}}(x)$ and $0 \not \in
\lambda_{\mathrm{ls}}(x)=(\limsup x) \upar$, which implies $\limsup
x =b>0$. By Lemma \ref{T1264} (b) and (c) we have $b_x=b$ and
$\bigvee_{y\prec x}\bigwedge_{z\prec y} \bigvee_{n \in \omega}
z_n=b$. Consequently, there exists $y \prec x$ and $c \in {\mathbb
B}^+$ such that $\bigwedge_{z\prec y} \bigvee_{n \in \omega} z_n=c$,
which implies
\begin{equation}\label{EQ1281d} \textstyle
\forall z \prec y~ \bigvee _{n \in \omega} z_n \geq c.
\end{equation}
\noindent{\it Claim 1.}
$\langle y_n\wedge c: n \in \omega \rangle$ is a  $\limsup$-stable  sequence.

\vspace{2mm}
\noindent {\it Proof of Claim 1.} First, by (\ref{EQ1281d}) and since $\langle y_n:n\geq k\rangle$ is a subsequence of $y$,
we have
$\limsup \; \langle y_n \wedge c: n\in \omega\rangle
= \bigwedge_{k \in  \omega} (\bigvee_{n\geq k} y_n) \wedge c
= \bigwedge_{k \in \omega}c
= c$.
Now we prove the same for an arbitrary subsequence
$\langle y_{f(k)} \wedge c: k\in \omega\rangle $ of $\langle y_n \wedge c: n\in \omega\rangle$, where
$f \in \omega^{\uparrow \omega}$. Clearly, $z=\langle y_{f(k)} : k\in \omega\rangle$ is a subsequence of $y$
and for each $l \in \omega$ we have
$\langle y_{f(k)} : k\geq l\rangle\prec y$, which, by (\ref{EQ1281d}), implies
$\bigvee_{k\geq l}y_{f(k)} \geq c$.
So, $\limsup \langle y_{f(k)}\wedge c :k \in \omega\rangle
=\bigwedge_{l\in \omega} \bigvee_{k\geq l}y_{f(k)}\wedge c
=\bigwedge_{l\in \omega}(\bigvee_{k\geq l}y_{f(k)})\wedge c
=\bigwedge_{l \in \omega}c=c$. Claim 1 is proved.

\vspace{2mm}
\noindent{\it Claim 2.} The set $M=\{n \in \omega: y_n \wedge c=0\}$ is finite.

\vspace{2mm}
\noindent{\it Proof of Claim 2.} Suppose that $M\in [\omega]^\omega$.
Then $\langle y_n \wedge c: n\in M\rangle$ is a subsequence of the sequence $\langle y_n \wedge c: n\in \omega\rangle $
and, clearly, $\limsup \langle y_n \wedge c: n\in M\rangle =0<c$,
which is impossible by Claim 1. Claim 2 is proved.

\vspace{2mm} By Claim 2, without loss of generality, we suppose that
$y_n\wedge c>0$, for each $n \in \omega$. By Theorem \ref{T1233a} we
have $\textstyle \overline{\{y_n \wedge c: n \in \omega\}} =c\upar
\cup \bigcup_{n\in \omega} (y_n \wedge c)\upar $ and this set is
closed in the space $\langle {\mathbb B},{\mathcal
O}_{\lambda_{\mathrm{ls}}} \rangle$, does not contain 0, but
contains each element of the sequence $\langle y_n\wedge c:n \in
\omega\rangle $. This implies $\textstyle 0 \not \in \Lim_{{\mathcal
O}_{\lambda_{\mathrm{ls}}}}\langle y_n \wedge c\rangle $.

On the other hand, since $y \prec x$ and $0 \in \Lim_{{\mathcal
O}_{\lambda_{\mathrm{ls}}}}(x)$, by (L2) we have $0 \in
\Lim_{{\mathcal O}_{\lambda_{\mathrm{ls}}}}(y)$. Since $y_n\wedge
c\leq y_n$, for each $n\in \omega $, by Lemma \ref{T4116}(c) we have
$\Lim_{{\mathcal O}_{\lambda_{\mathrm{ls}}}}(y) \subset
\Lim_{{\mathcal O}_{\lambda_{\mathrm{ls}}}}\langle y_n \wedge c
\rangle$ and, hence, $0 \in \Lim_{{\mathcal
O}_{\lambda_{\mathrm{ls}}}}\langle y_n \wedge c\rangle$. A
contradiction.
\kdok

\section{The algebras with $\lambda _{\mathrm{ls}}$ weakly topological}

By Theorem \ref{T1280}, if a complete Boolean algebra is not $(\omega ,2)$-distributive, the convergences
$\lambda_{\mathrm{ls}}$ and $\lambda_{\mathrm{li}}$ are not topological.
Now we show that they are weakly topological in algebras satisfying condition ($\hbar$).
The reader will notice that if in condition ($\hbar$) we replace ``$\limsup$" by ``$\liminf$", then we obtain
an equivalent condition, because $(\limsup x_n )'= \liminf x_n'$, for each sequence $x$ in ${\mathbb B}$.

\begin{te}\label{T1233}\rm
If ${\mathbb B}$ is a complete Boolean algebra satisfying condition
$(\hbar)$, then $\lambda _{\mathrm{ls} }$ and $\lambda _{\mathrm{li} }$ are weakly topological convergences.
\end{te}
\dok
We prove the statement for $\lambda _{\mathrm{ls} }$. The proof
for $\lambda _{\mathrm{li} }$ is dual. We show that for each
sequence $x$ in ${\mathbb B}$ and each $a\in {\mathbb B}$ we have
$a\in \lim_{{\mathcal O}_{\lambda _{\mathrm{ls}}}} x \Leftrightarrow
\forall y \prec x~ \exists z \prec y ~ \limsup z \leq a$. The
implication ``$\Leftarrow$" is Theorem \ref{T1261}(d). In order to
prove ``$\Rightarrow$" suppose that $a \in \lim_{{\mathcal
O}_{\lambda_{\mathrm{ls}}}} x$, $y \prec x$ and $\limsup z \not \leq
a$, for each subsequence $z \prec y$. By ($\hbar$), there is a
lim\,sup-stable sequence $z\prec y$. Then the set $K=\{ n\in \omega
: z_n \leq a \}$ is finite, since otherwise we would have $\limsup
\langle z_n: n \in K\rangle \leq a$. Thus w.l.o.g\ we can suppose
that $z_n \not \leq a$ for each $n \in \omega$. By Lemma
\ref{T1233a} we have
$$
\textstyle \overline{\{z_n:n \in \omega\}}=(\limsup z)\upar \cup
\bigcup_{n\in \omega} z_n\upar .
$$
Thus $a\in O={\mathbb B}\setminus \overline{\{z_n : n \in\omega\}} \in {\mathcal O}_{\lambda_{\mathrm{ls}}}$
and, since $O\cap \{z_n : n \in\omega\}=\emptyset$, we have
$a \not\in \lim_{{\mathcal O}_{\lambda_{\mathrm{ls}}}} z$. A contradiction, because $z\prec x$ and
$a \in \lim_{{\mathcal O}_{\lambda_{\mathrm{ls}}}} x$.
\kdok
\begin{ex}\rm\label{EX1204}
If ${\mathbb B}$ is a ccc complete Boolean algebra such that forcing
by ${\mathbb B}$ produces new reals, then, by Fact \ref{T1246} and
Theorems \ref{T1256}, \ref{T1280} and \ref{T1233}, the convergences
$\lambda_{\mathrm{s}},\lambda_{\mathrm{ls}}$ and
$\lambda_{\mathrm{li}}$ are weakly topological, but not topological.
In particular this holds for the Cohen algebra Borel$(2^\omega )/
{\mathcal M}$ and random algebra Borel$(2^\omega )/ {\mathcal Z}$,
where ${\mathcal M}$ and ${\mathcal Z}$ are the $\sigma$-ideals of
meager and measure-zero Borel sets, respectively.
\end{ex}
In the sequel, using the following lemma, we show that, on complete Boolean algebras belonging to a
large class, the convergence $\lambda_{{\mathrm{ls}}}$ is not weakly-topological.
\begin{lem}\rm \label{T1287}
Let ${\mathbb B}$ be a complete Boolean algebra, $x=\langle x_n:n\in \omega\rangle $
a sequence in ${\mathbb B}$ and $\tau_x=\{\langle \check{n},x_n\rangle : n\in \omega\}$ the corresponding ${\mathbb B}$-name
for a real. Then

(a) If $A$ is an infinite subset of $\omega$ and $f_A:\omega
\rightarrow A$ is the corresponding increasing bijection, then
 $ \| |\tau_x\cap \check{A}|=\check{\omega} \| =\limsup x\circ f_A$.

(b) The following conditions are equivalent:

(i)  $\forall f \in \omega^{\uparrow \omega} \; \exists g \in \omega^{\uparrow \omega} \; \limsup x\circ f \circ g =0$;

(ii) $\forall y\prec x \; \exists z \prec y \; \limsup z  =0$;

(iii) $\forall A \in [\omega]^{\omega}\; \exists B \in [A]^{\omega} \; \| |\tau_x \cap \check{B}| = \check{\omega }\| =0$.
\end{lem}
\dok
(a) Since $A=\{ f_A(n) : n\in \omega\}$ and $f_A$ is a bijection,
$\limsup x\circ f_A
=\bigwedge _{k\in \omega} \bigvee _{n\geq k} x_{f_A(n)}
= \| \forall k\in \check{\omega} \; \exists n\geq k \; f_A(n) \in \tau _x \|
= \| |\tau_x\cap \check{A}| =\check{\omega} \|$.

(b) The equivalence of (i) and (ii) is obvious.

(i) $\Rightarrow$ (iii) Let $A\in [\omega]^\omega$.  By (i), there is $g \in \omega^{\uparrow \omega}$
such that $\limsup x \circ f_A \circ g=0$.
Clearly $B=f_A[g[\omega]]\in [A]^\omega $ and $f_B=f_A \circ g$ so, by (a),
$\| |\tau_x \cap \check{B}| = \check{\omega }\| =\limsup x\circ f_A \circ g =0$.

(iii) $\Rightarrow$ (i) Let $f \in \omega^{\uparrow \omega}$ and
$A=f[\omega]$. By (iii),  there is $B \in [A]^\omega$ such that $\|
|\tau_x \cap \check{B}| = \check{\omega }\|=0$. Since $f^{-1}[B]\in
[\omega]^\omega$, there exists an increasing bijection $g:\omega
\rightarrow f^{-1}[B]$. From $B \subset f[\omega]$ it follows that
$f[g[\omega]]=f[f^{-1}[\omega]]=B$. So, by (a), $\limsup x\circ
f\circ g = \| |\tau _x \cap f[g[\omega ]]\check{\enspace}| =
\check{\omega }\| =0$ and (i) is proved.
\kdok
We remind the reader
that a set ${\mathcal T}\subset [\omega ]^{\omega }$ is called a
{\bf base matrix tree} iff $\langle {\mathcal T}, {}^* \!\! \supset
\rangle$ is a tree of height ${\mathfrak h}$ and ${\mathcal T}$ is a
dense set in the pre-order $\langle [\omega ]^{\omega }, \subset ^*
\rangle $. By a theorem of Balcar, Pelant and Simon (see
\cite{Balc0}), such a tree always exists. Clearly the levels of a
base matrix tree ${\mathcal T}$ are maximal almost disjoint families
and maximal chains in ${\mathcal T}$ are towers.
\begin{te}\rm \label{T1288}
If  ${\mathbb B}$ is a complete Boolean algebra satisfying
$1\Vdash_{\mathbb B} ({\mathfrak h}^V)^{\check{~}} <{\mathfrak t}$ and cc$({\mathbb B})>2^{\mathfrak h}$,
then $\lambda_{{\mathrm{ls}}}$ is not a weakly-topological convergence on ${\mathbb B}$.
\end{te}
\dok
Let  ${\mathcal T}$ be a base matrix tree and $\Br({\mathcal T})$ the set of all maximal branches of ${\mathcal T}$.
Since the levels of ${\mathcal T}$ are of size $\leq {\mathfrak c}$ and the height of ${\mathcal T}$ is
${\mathfrak h}$, for $\kappa =| \Br({\mathcal T})|$ we have $\kappa \leq {\mathfrak c}^{\mathfrak h}=2 ^{\mathfrak h}$
and we take an enumeration $\Br({\mathcal T})= \{ T_\alpha : \alpha < \kappa \}$.
Since $1\Vdash ({\mathfrak h}^V)^{\check{~}} <{\mathfrak t}$, for each $\alpha < \kappa $ we have
$1\Vdash |\check{T_\alpha }|<{\mathfrak t}$
and, consequently,
$1\Vdash \exists X \in [\check{\omega}]^{\check{\omega}} \; \forall B \in \check{T_\alpha }\; X \subset^*B$
so, by the Maximum Principle (see \cite[p.\ 226]{Kun}) there is
a name $\sigma_\alpha \in V^{\mathbb B}$ such that
\begin{equation}\label{EQ1290a}
1 \Vdash \sigma_\alpha  \in [\check{\omega}]^{\check{\omega}} \wedge
\forall B \in T_{\alpha } \;\; \sigma_\alpha \subset^*B.
\end{equation}
Let $\{b_\alpha: \alpha < \kappa \}$  be a maximal antichain in  ${\mathbb B}$.
By the Mixing lemma (see \cite[p.\ 226]{Kun})
there is a name $\tau \in V^{\mathbb B}$ such that
\begin{equation}\label{EQ1290b}
\forall \alpha <\kappa \;\; b_\alpha\Vdash \tau=\sigma_\alpha,
\end{equation}
and, clearly,  $1\Vdash \tau\in [\check{\omega}]^{\check{\omega}}$.
Let us define $x_n=\|\check{n}\in \tau\|$, $n \in \omega$.
Then for the corresponding name  $\tau _x = \{ \langle \check{n}, x_n \rangle : n\in \omega \}$ we have
\begin{equation}\label{EQ1290c}
1\Vdash \tau=\tau_x.
\end{equation}
\noindent
{\it Claim 1.} $0 \not \in \lambda^*_{{\mathrm{ls}}}(x)$.

\vspace{2mm} \noindent {\it Proof of Claim 1:} We prove that
$\neg\forall y \prec x ~\exists z \prec y~ \limsup z =0$ that is, by
Lemma \ref{T1287}(b), $\exists A \in [\omega]^\omega \; \forall B\in
[A]^\omega \; \| |\tau_x \cap \check{B}|=\check{\omega}\|>0$. In
fact, we show more:
\begin{equation}\label{EQ1290d}
\forall B \in [\omega]^\omega \;\; \| |\tau_x \cap \check{B}|=\check{\omega}\|>0.
\end{equation}
Let $B\in [\omega]^\omega$. Since ${\mathcal T}$ is a dense subset of $\langle [\omega]^\omega,\subset^*\rangle $
there is $C\in {\mathcal T}$ such that $C \subset^* B$. Let $T_\alpha$ be a branch
in ${\mathcal T}$ such that $C \in T_\alpha$. Then, by
(\ref{EQ1290b}) and (\ref{EQ1290c}) we have
$b_\alpha \Vdash \tau_x=\sigma_\alpha$, and by (\ref{EQ1290a})
$1 \Vdash \sigma_\alpha \subset^* C$, so $b_\alpha \leq \| |\tau_x \cap
\check{B}|=\check{\omega}\|$.

\vspace{2mm}
\noindent
{\it Claim 2.} $0  \in \Lim_{{\mathcal O}_{\lambda_{{\mathrm{ls}}}}}(x)$.

\vspace{2mm}
\noindent
{\it Proof of Claim 2:} On the contrary, suppose that there
are $F\in {\mathcal F}_{\lambda_{{\mathrm{ls}}}}$ and  $A\in [\omega]^\omega$ such that $0 \not \in F$
and $\{x_n:n \in A\} \subset F$.
Since ${\mathcal T}$ is dense in $\langle [\omega]^\omega,\subset^*\rangle $, there
is $C\in {\mathcal T}$ such that $C\subset^* A$ and, clearly, there is $\alpha <\kappa$ such that $C \in T_\alpha$.
$T_\alpha$ is a tower of type $\lambda \leq \mathfrak{h}$, so
$T_\alpha= \{B_\xi :\xi < \lambda \}$, where $B_\zeta \subsetneq^* B_\xi$, for $\xi < \zeta <\lambda$.
Let $C=B_{\xi_0}$ and, for $n\in \omega $, let
$$
D_n=B_{\xi_0+n}\setminus B_{\xi_0+n+1}.
$$
By Lemma \ref{T1287}(a), for each $n\in \omega $ we have
$\| |\tau_x\cap \check{D_n}|=\check{\omega}\| = \limsup x\circ f_{D_n}$.
Since $D_n \subset^* A$, almost all members of the sequence $x\circ f_{D_n}$ are elements of $F$ and, by Theorem \ref{T4105}(I),
$\| |\tau_x \cap \check{D_n} |=\check{\omega}\|\in F$.
So, by the same theorem,
$\limsup \| |\tau_x \cap \check{D_n} |=\check{\omega}\|\in F$.
Since
$\limsup \| |\tau_x \cap \check{D_n} |=\check{\omega}\|
= \| |\tau_x \cap \check{D_n} |=\check{\omega}\mbox { for infinitely many }n \in \omega\|$,
we will obtain a contradiction when we prove that
\begin{equation}\label{EQ1294}
\| |\tau_x \cap \check{D_n}| =\check{\omega}\mbox { for infinitely many }n \in \omega\|=0.
\end{equation}
Let $G$ be a ${\mathbb B}$-generic filter over $V$. Then there
exists $\beta<\kappa $ such that $b_\beta \in G$ and, by (\ref{EQ1290a}),(\ref{EQ1290b}) and (\ref{EQ1290c}),
\begin{equation}\label{EQ1295}
(\tau_x)_G \subset ^* B, \mbox{ for each } B\in T_\beta .
\end{equation}
First, if $\beta =\alpha$ then, by (\ref{EQ1295}), $|(\tau_x)_G \cap D_n|<\omega $,
for each $n\in \omega$.

\noindent
Second, if $\beta \neq \alpha$, we have two cases.

{\it Case 1:} $\exists E \in T_\beta \; \forall n\in \omega \; E \subset ^* B_{\xi _0 +n}$. Then $(\tau_x)_G \subset ^* E$
and for each $n\in \omega$ we have $|(\tau_x)_G \cap D_n|<\omega $.

{\it Case 2:} $\forall E \in T_\beta \; \exists n\in \omega \; E \not\subset ^* B_{\xi _0 +n}$.
Then, since ${\mathcal T}$ is a tree, there is the $\subset ^*$-maximum of the set $T_\beta \setminus T_\alpha$, say
$E'$ and, by the assumption, there is $n_0 \in \omega$ such that
$B_{\xi _0 +n_0 }\subset ^* E'$ or $|B_{\xi _0 +n_0}\cap E'|<\omega$.
Since $E' \not\in T_\alpha $,   $B_{\xi _0 +n_0 }\subset ^* E'$ is impossible, so $|B_{\xi _0 +n_0}\cap E'|<\omega$
and, hence, $|B_{\xi _0 +n}\cap E'|<\omega$, for each $n\geq n_0$. Since $(\tau _x)_G \subset ^* E'$ and
$D_n \subset  B_{\xi _0 +n}$, we have
$|(\tau_x)_G \cap D_n|<\omega $, for all $n\geq n_0$.

Thus $|(\tau_x)_G \cap D_n|<\omega $, for all but finitely many $n\in \omega$ and (\ref{EQ1294}) is true.
\kdok

The following example shows that there are very simple Boolean algebras such that the
question ``Is the convergence $\lambda_{\mathrm{ls}}$ on ${\mathbb B}$ weakly topological?" does not have an answer in ZFC.

\begin{ex}\rm\label{EX1205}
The statement ``The convergence $\lambda_{\mathrm{ls}}$ on the collapsing algebra ${\mathbb B}=\ro ({}^{<\omega }\omega _2)$
is weakly topological" is independent of ZFC. Since $\omega _2 ^{<\omega }=\omega _2$, the algebra ${\mathbb B}$ is
$\omega _3$-cc and collapses $\omega _2$ to $\omega$ in each generic extension.

If in the ground model $V$ we have $2^\omega =\omega _1$ and $2^{\omega _1}=\omega _2$ (in particular, if $V\models$ GCH)
then in $V$ we have ${\mathfrak h}=\omega _1$, cc$({\mathbb B})=\omega _3 > \omega _2 =2^{{\mathfrak h}}$ and
$1\Vdash _{\mathbb B} |({\mathfrak h}^V)^{\check{~}}|=\check{\omega }$. Thus, by Theorem \ref{T1288}, the convergence
$\lambda_{\mathrm{ls}}$ on ${\mathbb B}$ is not weakly topological.

On the other hand, if in $V$ we have ${\mathfrak t}\geq \omega _3$
(in particular, if $V\models {\mathrm{MA}} + {\mathfrak c}\geq \omega _3)$, then
${\mathbb B}$ is ${\mathfrak t}$-cc and, hence, satisfies condition $(\hbar )$ which, by Theorem \ref{T1233},
implies that the convergence
$\lambda_{\mathrm{ls}}$ on ${\mathbb B}$ is weakly topological.
\end{ex}


\footnotesize

\begin{thebibliography}{aa}
\bibitem{BGJ98}
      B.\ Balcar, W.\ Gl\'{o}wczy\'{n}ski, T.\ Jech,
      The sequential topology on complete Boolean algebras,
      Fund.\ Math.\ 155 (1998) 59--78.
\bibitem{BJP05}
      B.\ Balcar, T.\ Jech, T.\ Paz\'{a}k,
      Complete ccc Boolean algebras, the order sequential topology and a problem of von Neumann,
      Bull.\ Lond.\ Math.\ Soc.\ 37,6 (2005) 885--898.
\bibitem{BJ06}
      B.\ Balcar, T.\ Jech,
      Weak distributivity, a problem of von Neumann and the mistery of measurability,
      Bull.\ Symbolic Logic, 12,2  (2006) 241--266.
\bibitem{Balc0}
      B.\ Balcar, J.\ Pelant, P.\ Simon,
      The space of ultrafilters on N covered by nowhere dense sets,
      Fund.\ Math., 110 (1980) 11--24.
\bibitem{Douw84}
      E.K.\ van Douwen,
      The integers and topology,
      in: K.\ Kunen and J.E.\ Vaughan eds., Handbook of Set-theoretic Topology,
      North-Holland, Amsterdam, 1984, 111--167.
\bibitem{Eng85}
      R.\ Engelking,
      General Topology,
      P.W.N.\ Warszawa, 1985.
\bibitem{Far04}
      I.\ Farah,
      Examples of $\varepsilon$-exhaustive pathological submeasures,
      Fund.\ Math.\ 181 (2004) 257--272.
\bibitem{Jec97}
      T.\ Jech,
      Set Theory, 2.\ corr.\ ed.,
      Springer, Berlin, 1997.
\bibitem{Kun}
      K.\ Kunen,
      Set Theory,
      An Introduction to Independence Proofs,
      (North-Holland, Amsterdam, 1980).
\bibitem{KuPa07}
      M.S.\ Kurili\'{c}, A.\ Pavlovi\'{c},
      A posteriori  convergence  in complete Boolean algebras with the sequential topology,
      Ann.\ Pure Appl.\ Logic  148,1-3 (2007) 49--62.
\bibitem{KuPaNSJOM}
      M.S.\ Kurili\'{c}, A.\ Pavlovi\'{c},
      Some forcing related convergence structures on complete Boolean algebras,
      Novi Sad J.\ Math.\ 40,2 (2010) 77--94.
\bibitem{KuPaDEBR}
      M.S.\ Kurili\'{c}, A.\ Pavlovi\'{c},
      The convergence of the sequences coding the ground model reals,
      submitted.
\bibitem{KuTo}
      M.S.\ Kurili\'c, S.\ Todor\v cevi\'c,
      Property ($\hbar$) and  cellularity of complete Boolean algebras,
      Arch.\ Math.\ Logic, 48,8 (2009) 705--718.
\bibitem{Mah47}
      D.\ Maharam,
      An algebraic characterization of measure algebras,
      Ann.\ of Math., 48 (1947) 154--167.
\bibitem{Sco81}
      R.D.\ Mauldin (ed.),
      The Scottish Book (Mathematics from the Scottish Caf\'{e}),
      Birkh\"auser, Boston MA, 1981.
\bibitem{Tal1}
      M.\ Talagrand,
      Maharam's problem,
      C.\ R.\ Acad.\ Sci.\ Paris, Ser.\ I, 342 (2006) 501--503.
\bibitem{Tal2}
      M.\ Talagrand,
      Maharam's problem,
      Ann.\ of Math., 168,3  (2008) 981--1009.
\bibitem{Tod04}
      S.\ Todor\v{c}evi\'{c},
      A problem of von Neumann and Maharam about algebras supporting continuous submeasures,
      Fund.\ Math., 183,2 (2004) 169--183.
\bibitem{Vel05}
      B.\ Veli\v{c}kovi\'{c},
      ccc forcing and splitting reals,
      Israel J.\ Math., 147 (2005) 209--220.
\end{thebibliography}
\end{document}